\documentclass[12pt, reqno]{amsart}
\usepackage[arrow,matrix,curve]{xy}
\usepackage{graphicx}

\newtheorem{theorem}{Theorem}[section]
\newtheorem{proposition}[theorem]{Proposition}
\newtheorem{lemma}[theorem]{Lemma}

\newtheorem{sublemma}[theorem]{Sublemma}

\theoremstyle{definition}

\newtheorem{remark}[theorem]{Remark}

\newcommand{\R}{\mathbb R} 
\newcommand{\Z}{\mathbb Z}

\numberwithin{equation}{section}
\numberwithin{figure}{section}

\begin{document}

\author[C.~Croke]{Christopher Croke$^+$} \address{ Department ofbb:
Mathematics, University of Pennsylvania, Philadelphia, PA 19104-6395
USA} \email{ccroke@math.upenn.edu} \thanks{Supported in part by NSF grant DMS 10-03679 and an Eisenbud Professorship at M.S.R.I.}

\title[Scattering rigidity] {Scattering rigidity with trapped geodesics
}


\begin{abstract}
We prove that the flat product metric on $D^n\times S^1$ is scattering rigid where $D^n$ is the unit ball in $\R^n$ and $n\geq 2$.

The scattering data (loosely speaking) of a Riemannian manifold with boundary is map $S:U^+\partial M\to U^-\partial M$ from unit vectors $V$ at the boundary that point inward to unit vectors at the boundary that point outwards.  The map (where defined) takes $V$ to $\gamma'_V(T_0)$ where $\gamma_V$ is the unit speed geodesic determined by $V$ and $T_0$ is the first positive value of $t$ (when it exists) such that $\gamma_V(t)$ again lies in the boundary.

We show that any other Riemannian manifold $(M,\partial M,g)$ with boundary $\partial M$ isometric to $\partial(D^n\times S^1)$ and with the same scattering data must be isometric to $D^n\times S^1$.

This is the first scattering rigidity result for a manifold that has a trapped geodesic.  The main issue is to show that the unit vectors tangent to trapped geodesics in $(M,\partial M,g)$ have measure 0 in the unit tangent bundle.

\end{abstract}

\maketitle

\section{Introduction}
In this paper we prove scattering rigidity (see below) for a number of compact Riemannian manifolds with boundary that have trapped geodesics.  A geodesic $\gamma(t)$ in a manifold with boundary is trapped if its domain of definition is unbounded.

Consider a compact Riemannian manifold $(M,\partial M,g)$ with boundary $\partial M$ and metric $g$.  We will let $U^+\partial M$ represent the space of inwardly pointing unit vectors at the boundary.  That is $V\in U^+\partial M$ means that $V$ is a unit vector based at a boundary point and $\langle V,\eta^+\rangle \geq 0$ where $\eta^+$ is the unit vector of $M$ normal to $\partial M$ and pointing inward.  Similarly we let $U^-\partial M$ represent the outward vectors.  Note that $U^+\partial M\cap U^-\partial M=U(\partial M)$ the unit tangent bundle of $\partial M$.

For $V\in U^+\partial M$ let $\gamma_V(t)$ be the geodesic with $\gamma'(0)=V$.  We let $TT(V)\in [0,\infty]$ (the travel time) be the first time $t>0$ when $\gamma_V(t)$ hits the boundary again.  If $\gamma_V(t)$ never hits the boundary again then $TT(V)=\infty$ (i.e $\gamma_V$ is trapped) while if either $\gamma_V(t)$ does not exist for any $t>0$ or there are arbitrarily small values of $t>0$ such that $\gamma(t)\in \partial M$ then we let $TT(V)=0$.  Note that $TT(V)=0$ implies that $V\in U(\partial M)$.

The scattering map ${\emph{S}}:U^+\partial M\to U^-\partial M$ takes a vector $V\in U^+\partial M$ to the vector $\gamma'(TT(V))\in U^-\partial M$.  It will not be defined when $TT(V)=\infty$  and will be $V$ itself when $TT(V)=0$.
If another manifold $(M_1,\partial M_1,g_1)$ has isometric boundary to $(M,\partial M,g)$ in the sense that $(\partial M,g)$ (g restricted to $\partial M$) is isometric to $(\partial M_1,g_1)$ then we can identify $U^+\partial M_1$ with $U^+\partial M$ and $U^-\partial M_1$ with $U^-\partial M$.
We say that $(M,\partial M,g)$ and $(M_1,\partial M_1,g_1)$ have the same scattering data if they have isometric boundaries and under the identifications given by the isometry they have the same scattering map.  If in addition the travel times $TT(V)$ coincide then they are said to have the same lens data.

A compact manifold $(M,\partial M,g)$ is said to be scattering (resp. lens) rigid if for any other manifold $(M_1,\partial M_1,g_1)$ with the same scattering (resp. lens) data there is an isometry from $M_1$ to $M$ that agrees with the given isometry of the boundaries.

\begin{theorem}
\label{3+dims}
For any $n\geq 2$ the flat product metric on $D^n\times S^1$ is scattering rigid where $D^n$ is a ball in $\R^n$.
\end{theorem}

\begin{figure}
\label{cylinders}
\includegraphics[width=50mm]{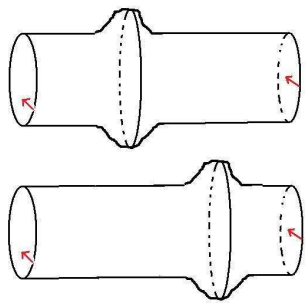}
\caption{not isometric but same scattering and lens data}

\end{figure}

The fact that not all manifolds are scattering rigid was pointed out in \cite{Cr91}.  For $\frac 1 4 >\epsilon>0$ let $h(t)$ be a small smooth bump function which is 0 outside $(-\epsilon ,\epsilon)$ and positive in $(-\epsilon,\epsilon)$.   For $s\in (-1+2\epsilon,1-2\epsilon)$ consider surfaces of revolution $g_s$ with smooth generating functions $F_s(t)=1+h(s+t)$ for $t\in [-1,1]$.  These surfaces of revolution look like flat cylinders with bumps on them that are shifted depending on $s$ but otherwise look the same (see figure \ref{cylinders}).  The Clairaut relations show that, independent of $s$, geodesics entering one side with a given initial condition exit out the other side after the same distance at the same point with the same angle.  Hence all metrics have the same scattering data (and in fact lens data) but are not isometric.  A much larger class of examples was given in section 6 of \cite{Cr-Kl94}.  All of the examples have in common that there are trapped geodesics.

The scattering and lens rigidity problems are closely related to other inverse problems.  In particular the boundary rigidity problem is equivalent to the lens rigidity question in the Simple and SGM cases.  See \cite{Cr91} and \cite{Cr04} for definitions and relations to some other problems.  There is a vast literature on these problems (see for example \cite{Be83,Bu-Iv06,Cr91,Cr90,Gr83,Mi81,Mu77,Ot90,Pe-Sh88,Pe-Ul05,St-Uh09}).  Most of the results in these papers concern manifolds with no trapped geodesics.  An exception is \cite{St-Uh09} where they prove local scattering rigidity (i.e. if the two metrics are in a particular $C^k$ neighborhood they must be isometric) for a class of Riemannian manifolds that includes those discussed in this paper.  However, to date all of the global rigidity results concern manifolds without trapped geodesics.  The results in this paper constitute the first examples of (global) scattering rigid manifolds that have trapped geodesics.

The key difficulty in our case is to show that the set of unit vectors tangent to trapped geodesic rays in the metric $g_1$ has measure 0 in the unit tangent bundle.  This allows us (with an application of Santal\'o's formula) to conclude that $g$ and $g_1$ have the same volumes.  Since the metric $g$ has a real factor (i.e. $D^n\subset \R\times \R^{n-1}$) we can use a result from \cite{Cr-Kl98} to complete the argument.  In fact, the argument in Theorem \ref{3+dims} extends (see section \ref{gens}) to the case where $D^n$ above is replaced by a ball in $\R\times N^{n-1}$ where $N$ is a complete simply connected Riemannian manifold with nonpositive curvature.  (In fact with more work one could extend this to the case of no conjugate points but we chose not to give the slightly different arguments here.)

One case that was not dealt with in Theorem \ref{3+dims} is the two dimensional case, namely the flat cylinder $[-1,1]\times S^1$ and the M\"obius strip.  There are ways in which this case is easier and ways in which it is harder.  The major differences are that the scattering data does not determine the lens data and we cannot conclude that the $C^\infty$ jets of the metrics agree at the boundary.  The problem of lens rigidity in the two dimensional case will be taken up in a future paper with Pilar Herreros.  In particular, it turns out that the M\"obius strip is not scattering rigid if $(M_1,\partial M_1,g_1)$ is allowed to be $C^1$.

The author would like to thank Gunther Uhlmann who first posed the problem of the rigidity of $D^2 \times S^1$ to him some years ago, to Haomin Wen for pointing out the Eaton Lens example, and Pilar Herreros for a careful reading of earlier drafts.

\section{The $D^n\times S^1$ for $n\geq 2$ case}
\label{highdim}

In this section we prove Theorem \ref{3+dims}.  We consider generalizations in Section \ref{gens}.  Throughout the section $n\geq 2$ and $g$ will be the standard flat product metric on $M=D^n\times S^1$.  For concreteness we will take $D^n$ to be the unit ball and $S^1$ to have length $2\pi$.  $(M_1,S^{n-1}\times S^1,g_1)$ will be another Riemannian metric on a manifold $M_1$ whose boundary is isometric to that of $M$.  We use this isometry to identify the two boundaries.  We assume that $g_1$ has the same scattering data as $g$.

We do not a-priori assume that $M_1$ is diffeomorphic to $M$.  For each $p\in \partial M=\partial M_1$ we let $\tau_p\subset \partial M=\partial M_1$ be the closed curve in vertical (i.e. $S^1$) direction.

\begin{lemma}
\label{lensdata}
$g_1$ has the same lens data as $g$.
\end{lemma}

{\bf Proof:}  This is an application of the first variation formula.  For $V\in U^+\partial M$ let $G(V)= L(\gamma_{1 V})-L(\gamma_{V})$.  We need to show that $G(V)=0$ for all $V$.  A smooth curve of initial conditions $s\mapsto V(s)$ in the interior of $U^+\partial M$ gives rise to smooth variations $\gamma_{V(s)}$ through unit speed geodesics in $M$ and $\gamma_{1V(s)}$ through unit speed geodesics in $M_1$ whose initial and final tangents agree.  (Note that this uses the convexity of the boundary since for more general manifolds with boundary there may be a discontinuous jump in the endpoints of geodesics.) The first variation formula (along with the fact that the metrics agree at the boundary) tells us that $\frac d{ds} L(\gamma_{V(s)})= \frac d{ds} L(\gamma_{1V(s)})$ .  Hence $G(V(s))= L(\gamma_{1V(s)})-L(\gamma_{V(s)})$ is independent of $s$.  Since $U^+\partial M$ is connected, G is a constant $C$.  Further when $V$ approaches a non vertical vector (i.e. one not tangent to the $S^1$ factor) in $\partial(U^+\partial M)=U\partial M$ then $L(\gamma_V)$ approaches $0$ and hence $C>0$.  If we knew $L(\gamma_{1 V})$ approached $0$ for this or any sequence then $C=0$ and the lemma would follow.  We now show that this must happen.

If this is not the case (i.e. if $C>0$) then the boundary must be {\em concave} everywhere and further when $V_i$ approaches a non vertical vector $V$ in $\partial(U^+\partial M)=U\partial M$ as above then $L(\gamma_{1 V})$ approaches $C$ and the limiting geodesic $\gamma_{1V}$ must be a closed geodesic of length $C$ with initial (and final) tangent vector $V$.  By taking limits this is also true for the vertical vectors.  (This certainly looks unlikely to happen.  However, there is an example \cite{H-H-L06} - an Eaton lens - of a manifold with a singularity having the same scattering data as the flat 2-disc but different lens data.  In that case the boundary is a closed geodesic.)

Note that in our case (since the dimension of the boundary is at least 2) all of these closed geodesics $\gamma_{1V}$ for $V$ based at a boundary point $p$ are homotopic to each other since there is a curve of initial tangent vectors $V_t$ between any two initial tangents (i.e. the homotopy is through the curves $\gamma_{1V_t}$). In particular, they are all homotopic to their negatives (i.e. running around in the opposite direction).  Thus going twice around such a geodesic is a contractible curve.

In fact, the closed geodesic $\gamma_V$ of length $C$ tangent to vertical direction $V$ at $p\in \partial M$ is a multiple of $\tau_p$.  To see this let $x=\tau_p(t)$ be a point on $\tau_p$ close to $p$ (say $t<\frac \pi 2$) and $x_i$ a sequence of points on $\partial M$ approaching $x$ but not on $\tau_p$.  Let $\gamma_{V_i}$ be the minimal $g$ geodesics from $p$ to $x_i$ with $x_i=\gamma_{V_i}(t_i)$.  We see that $V_i$ approaches $V$ and $t_i$ approaches $t$.  Looking at the other metric we see $x_i=\gamma_{1 V_i}(t_i+C)$ and hence by taking limits $x=\gamma_{1 V}(t+C)$.  Since this is true for all $p$ and $x=\tau_p(t)$ for $t<\frac \pi 2$  $\frac \pi 2$ we see that $\tau_p(t)= \gamma_{1 V_i}(t+C)$ and hence $\tau_p$ is the $g_1$ geodesic with initial tangent $V$.  Thus $\gamma_{1 V_i}$ simply goes around $\tau_p$ a number of times (i.e. $C$ is an integer multiple of $2\pi$).  Thus we know that going twice as many times around $\tau_p$ yields a contractible curve.  Thus the next sublemma \ref{noncontractible} gives the desired contradiction.
\qed

\begin{sublemma}
\label{noncontractible}
No multiple of $\tau_p$ is contractible in $M_1$.
\end{sublemma}

{\bf Proof:} This argument is an oriented intersection number argument.  We begin by seeing that when $C>0$ then $M_1$ is orientable.  Fix $x$ on the boundary.  Then for every element of $\pi_1$ the shortest loop representing this class is either a geodesic loop in $M_1$ or partly runs along the boundary.  By the scattering data assumption the only geodesic loops at $p$ are those that start tangent to the boundary (i.e. the closed geodesics we are discussing).  All of these are homotopic to a multiple of $\tau_p$.  If the minimizing path runs along the boundary some of the time then (since when it leaves the boundary it must be tangent) the only parts not on the boundary are the closed geodesic loops again.  Hence every element of $\pi_1$ has a representative that lies in the boundary.  Since running along such curves does not change the orientation, $M_1$ must be orientable.

Chose an orientation reversing diffeomorphism $F:S^{n-1}\to S^{n-1}$.  This induces a map $H:S^{n-1}\times [0,1] \to M_1$ by
$$H(x,0)=(x,0)\in \partial M_1=S^{n-1}\times S^1\ \ \ \ \ \ \ H(x,1)=(F(x),\frac \pi {10}).$$
For fixed $x$, $H(x,t)$ is the $g_1$ geodesic (parameterized proportional to arclength) with the same initial tangent as the minimizing $g$ geodesic from $(x,0)$ to $(F(x),\frac \pi {10})$ (and hence has length $C$ longer).  Note that if $F(x)=x$ then the geodesic $H(x,t)$ may wrap around $\tau_{(x,0)}$ many times before ending at $(x,\frac \pi {10})$.)

Although $H$ may not be transverse to the boundary that is easy to fix.  For some small $\epsilon>0$ we can parameterize the $\epsilon$ neighborhood $N_\epsilon$  in $M_1$ of $\partial M_1$ as $[0,\epsilon)\times S^{n-1}\times S^1$ (where $\partial M_1=\{0\}\times S^{n-1}\times S^1$).  Further we have a diffeomorphism $D_\epsilon:M_1\to M_1-N_\epsilon$ such that $D_\epsilon(0,x,\theta)=(\epsilon,x,\theta)$.  Thus we can take $\tilde H:S^{n-1}\times [-\epsilon,1+\epsilon]\to M_1$ by letting $\tilde H (x,t)= (t,x,0)$ for $t\leq 0$, $\tilde H(x,t)=D_\epsilon(H(x,t))$ for $0\leq t \leq 1$, and let $\tilde H(x,t)=(1+\epsilon-t,F(x),\frac \pi {10})$ for $t\geq 1$.  This is now transverse to $\partial M_1$ and we can tweak this to make it smooth.  We can double this picture in the manifold double $M_1\times M_1$ to get a smooth map $\bar H: S^{n-1}\times S^1\to M_1\cup M_1$ which is transverse to the curve $\tau_p$.  $\tau_p$ passes through each of $S^{n-1}\times \{0\}\in \partial M_1$ and $S^{n-1}\times \{\frac \pi {10}\}\in \partial M_1$ once each time around and our choice of orientations guarantees that the two contributions to the intersection number of $\bar H$ with $\tau_p$ have the same sign.  Thus $\bar H$ has a nonzero intersection with any multiple of $\tau_p$.  Thus the homology class of any multiple of $\tau_p$ is non zero in $M_1\cup M_1$ and hence in $M_1$.

\qed

Since the lens data and hence distances between nearby boundary points agree, the $C^\infty$ jets of $g$ and $g_1$ also must agree at the boundary.  This follows from \cite{L-S-U03},\cite{Uh-Wa03}, or \cite{Zh11} since, for the flat metric $g$, the second fundamental form of the boundary has a positive eigenvalue at every point.  (Note that this argument wont work in the two dimensional case $n=1$.)  This in particular means that we can glue $(\R^n-D^n)\times S^1$ along the boundary of $M_1$ to yield a smooth metric $M_1^{ext}$ which is isometric to $\R^n\times S^1$ outside of $M_1$.

\begin{lemma}
\label{fundgroup}
$\pi_1(M_1)=\Z$ and the generator is represented by the $S^1$ factor of the boundary.
\end{lemma}

{\bf Proof:} To see this fix a base point $p\in \partial M_1$.  The lens data being the same tells us that the only geodesic loops at a $p$ are the multiples of $\tau_p$.  Further the convexity of the boundary guarantees that there is at least one geodesic loop in each homotopy class (the shortest curve in that class).  In particular we see as before that $M_1$ is orientable.  The lemma now follows from the proof of sublemma \ref{noncontractible}.

\qed

\bigskip

Thus the Riemannian universal cover $\widetilde{M_1}$ of $M_1$ sits naturally in $\widetilde{M_1^{ext}}$ the universal cover of $M_1^{ext}$ and further $\widetilde{M_1^{ext}} - \widetilde{M_1}$ is isometric to $(R^n-D^n)\times \R$.  Also $\widetilde M = D^n\times \R$ has the same scattering data as $\widetilde M_1$.  We will slightly abuse notation and call the metrics on the universal covers $g$ and $g_1$ as well.

\begin{lemma}
\label{geosminimize}
An $M_1$ geodesic $\gamma$ between boundary points is the shortest path in its homotopy class (rel boundary points).
\end{lemma}

{\bf Proof:}  This is the same as saying that such geodesics in the universal cover are the minimizing paths between the endpoints.  This is true for $\widetilde M$ where there is a unique geodesic between any two boundary points.  Thus there is also a unique geodesic between boundary points in $\widetilde M_1$ which must thus be the minimizing geodesic.

\qed

\bigskip

In fact, this implies that all geodesic segments in $\widetilde{M_1^{ext}}$ are minimizing except possibly in the case that they are segments of geodesics trapped in $\widetilde{M_1}$.  If $p$ and $q$ are points in $\widetilde{M_1^{ext}}-\widetilde{M_1}=\R^{n+1}-D^n\times \R$ then this implies that $d_1(p,q)=d(p,q)$.  In particular, for $p\in \widetilde{M_1^{ext}}-\widetilde{M_1}$ all geodesics from $p$ minimize.  Hence the exponential map is a diffeomorphism which allows us to conclude that not only are the fundamental groups the same but $M_1$ is diffeomorphic to $M$.

A geodesic $\gamma_{1V}$ will either be trapped or coincide with an oriented Euclidean line $L_V$ outside $\widetilde{M_1}$.  By the direction of $L_V$ we mean the oriented line through the origin parallel to $L_V$.  $L_V$ will be called ``positive'' if it makes a (strictly) positive inner product with the upward vertical direction.

There are two cases that are exceptional.  These are trapped geodesics and vertical geodesics (i.e. $\{x\}\times \R$ for $x\in R^n-D^n$).  We will exclude both these cases by the phrase ``$L_v$ is not vertical''.  For $p\in \widetilde{M_1^{ext}}$ we let ${\bf A}(p)=\{V\in U_p\widetilde{M_1^{ext}}|L_V \textrm{ is not vertical}\}$.  Note that for $p\in \widetilde{M_1^{ext}}-\widetilde{M_1}$ we have ${\bf A} (p)$ is just the unit sphere with the north and south pole removed.  {$\bf A$ will represent the union of the ${\bf A}(p)$.

\begin{lemma}
\label{vertical}
If $V_i\in {\bf A}(p)$ converges to a vector $V\in U_p-{\bf A}(p)$ then the directions of the lines $L_{V_i}$ become vertical.
\end{lemma}

{\bf Proof:} Assume this is not the case.  Then there is a subsequence of the $V_i$ (which we will still call $V_i$) such that the directions of the lines $L_{V_i}$ converge to a non vertical direction $L$.  We claim that a subsequence of these $L_{V_i}$ converge to a line $L_W$.  To see this we only need to note that the lines $L_{V_i}$ intersect a common compact set.  Now the length of $\gamma_{V_i}\cap\widetilde M_1$ is the same as the length of $L_{V_i}\cap D^n\times \R\subset \R^{n+1}$ which is uniformly bounded above (say by $B$) since the directions of the $L_{V_i}$ converge to $L$ which is not vertical.  Thus all the $L_{V_i}$ intersect the boundary of $\widetilde M_1$ inside the compact ball about $p$ of radius $B$.  This means that the geodesics $\gamma_{V_i}$ (for the subsequence) converge to $L_V$ outside $\widetilde M_1$ but they converge to $\gamma_V$ which is supposed to be trapped.  This yields the desired contradiction.

\qed

We next see that even though $\widetilde{M_1^{ext}}$ might a-priori have conjugate points (along geodesics trapped in $\widetilde{M_1}$), Busemann functions along rays where $L_V$ is not vertical behave like those in manifolds without conjugate points.  In particular they are $C^{1,1}$ smooth, $|\nabla b_{1V}|=1$, and the Lipshitz constant for $\nabla b_{1V}$ is uniformly bounded.

For $V\in U\widetilde{M_1^{ext}}$ such that $\gamma_V(t)$ minimizes for all positive $t$, let $b_{1V}:\widetilde{M_1^{ext}}\to \R$ be the Busemann function defined by $V$, i.e.
$$b_{1V}(p)=lim_{t\to\infty} d_1(p,\gamma_v(t))-t.$$

Since $d_1(p,q)=d(p,q)$ when $p$ and $q$ are points in $\widetilde{M_1^{ext}}-\widetilde{M_1}$, $b_{1V}$ coincides with the Euclidean $b_V$ outside $\widetilde{M_1}$ as long as $\gamma_{1V}$ is not trapped.  That is $b_{1V}$ will coincide with the height function (up to a constant) in the direction $L_V$.

For all reals $s$ we will let $H_V(s)=\{p\in \widetilde{M_1^{ext}}| b_{1V}(p)=s \}$ be the $s$ level set of $b_{1V}$. Of course, outside of $\widetilde{M_1}$, $H_V(s)$ is just a hyperplane perpendicular to $L_V$.

\begin{lemma}
\label{busemannsmooth}
For all $V\in {\bf A}$, $b_{1V}$ is $C^{1,1}$ and the Lipshitz constant of $\nabla b_{1V}$ is bounded by a constant independent of $V$.

\end{lemma}

{\bf Proof:}  The proof is the usual proof that such a statement holds on manifolds without conjugate points.  This is done by showing that the approximating functions, $f_t(p)=d_1(p,\gamma_{1V}(t))-t$, are $C^\infty$ have $|\nabla f_t|=1$ and have uniformly bounded Hessian.  If $\gamma_{1V}(t)\in \widetilde{M_1^{ext}}-\widetilde{M_1}$ then all geodesics from $\gamma_{1V}(t)$ minimize so the distance function from $\gamma_{1V}(t)$ is $C^\infty$ for large $t$.  The fact that $|\nabla f_t|=1$ is clear.  The uniform control on the Hessian is also the same as we will see.  Fix a number $r$ less than the convexity radius of ${M_1^{ext}}$ - which exists since $M_1^{ext}-M_1$ is flat.  Since there is a compact set $K\subset M_1^{ext}$ of base points such that for $q\notin K$ the ball $B(q,r)$ is flat we conclude that the eigenvalues of the second fundamental forms of the boundaries of $B(q,r)$ are uniformly bounded above and below independent of $q$.  This same bound applies to balls in the universal cover $\widetilde{M_1^{ext}}$. Now to bound the Hessian of $f_t$ at $q\in \widetilde{M_1^{ext}}$ let $\tau(s)$ be the (unique) geodesic from $\gamma_{1V}(t)$ to $q=\tau(s_0)$ (we can assume $s_0\gg r$ since we will be taking the limit as $t\to \infty$).  Then by the triangle inequality the level set of $f_t$ at $q$ (i.e. $\partial B(\gamma_{1V}(t),s_0)$) lies outside both $B(\tau(s_0-r),r)$ and $B(\tau(s_0+r),r)$ which are tangent to the level set at $q$.  Hence the second fundamental forms of the level sets are uniformly bounded and hence so is the Hessian.  Thus the lemma follows.

\qed

The usual properties of Busemann functions (see \cite{Es77} for basic properties of Busemann functions) tell us that if $W(p)=\nabla b_{1V}(p)$ then $\gamma_{1W}'(t) = \nabla b_{1V}(\gamma_{1W}(t))$ for all $t$.  Hence if $\gamma_{1W}$ is not trapped then $L_W$ will be parallel to $L_V$.  A straightforward open and closed argument shows that for all $p$, $\gamma_{1W(p)}$ is not trapped.

\begin{lemma}
\label{levelsets}
Let $V$ and $W$ in ${\bf A}$ be such that $L_V$ and $L_W$ are positive and not parallel to each other.  Then for any given $s$ the maximum and minimum values of $b_{1V}$ on the compact $H_W(s)\cap \widetilde{M_1}$ are achieved on the boundary of $H_W(s)\cap \widetilde{M_1}$.

\end{lemma}

{\bf Proof:}  We first note that $H_W(s)\cap \widetilde{M_1}$ is compact.  Indeed, if $D$ is the diameter of $M_1$ then for every $p\in H_W(s)\cap \widetilde{M_1}$ there is a point $q\in \partial(\widetilde{M_1})$ with $d(p,q)\leq D$.  Since $b_{1W}$ has Lipshitz constant 1, we know that $s-D \leq b_{1w}(q)\leq s+D$ and hence $p$ lies in the (compact) set of points that are at distance $\leq D$ from the compact (since $W$ is not horizontal since it is positive) set of boundary points $\{q\in \partial \widetilde{M_1}| s-D\leq b_{1W}(q)\leq s+D\}=\{q\in \partial (B^{n-1}\times \R)| s-D\leq b_{W}(q)\leq s+D\}$.

If the maximum (or minimum) value of $b_{1V}$ occurs in the interior then $\nabla b_{1V}$ must be perpendicular to $H_W(s)$ there and hence coincides with $\pm \nabla b_{1W}$ at that point.    However this contradicts the condition that $L_V$ and $L_W$ are positive and not parallel.

\qed

We now see that if $V$ and $W$ in ${\bf A}$ are such that $L_V$ and $L_W$ are parallel then $b_{1V}-b_{1W}$ is constant.  Since they agree with the height functions outside $\widetilde{M_1}$, $b_{1V}-b_{1W}=C$ outside $\widetilde{M_1}$.  But then they must also differ by $C$ along any geodesic whose corresponding line is parallel to $L_V$ and $L_W$.  But since such a geodesic passes though every point $p\in \widetilde{M_1}$ (i.e. take the geodesic in the direction of $\nabla b_{1V}(p)$) this says $b_{1V}-b_{1W}=C$ everywhere.  In particular for every $p$ there is a unique geodesic passing through $p$ and parallel to a given line.  This gives a natural identification of ${\bf A}(p)$ with the space of non vertical directions.

\begin{proposition}
\label{mainhigher}
Through every $p\in \widetilde M_1$ there is exactly one trapped geodesic.
\end{proposition}

{\bf Proof:} This is equivalent to showing that for every $p\in \widetilde{M_1^{ext}}$, ${\bf A}(p)$ consists of the unit sphere $U_p\widetilde{M_1^{ext}}$ minus a pair of antipodal points.  Assume that $p$ is a point with more trapped geodesics.  Note that if there is a trapped half geodesic at $p$ then the other half must also be trapped by the assumption that the scattering data coincides with the flat case.  Thus the tangent directions to trapped geodesics come in antipodal pairs.  Of course there is at least one trapped geodesic through $p$ since one could take the limit of a subsequence of geodesics from $p$ to a boundary points $q_i$ where $q_i$ runs off to infinity.  We only need to consider $p$ in the interior of $\widetilde{M_1}$.

A limiting half geodesic of a sequence of trapped half geodesics starting at $p$ will be a half geodesic starting at $p$ that stays in $\widetilde M_1$.  (In fact it stays in the interior since the only half geodesics in $\widetilde M_1$ that are tangent to the boundary are the vertical ones hence stay in the boundary.)   Thus the set of tangent directions to trapped geodesics (i.e. $U_p-{\bf A}(p)$) is closed in $U_p$ and thus ${\bf A}(p)$ is open and nonempty (by the correspondence with non vertical directions).  The set of boundary points of $U_p-{\bf A}(p)$ is thus non empty and if it consisted of a single antipodal pair then $U_p-{\bf A}(p)$ would also be a single antipodal pair.  Thus there is a pair of distinct unit vectors $V$ and $W$ in $U_p-{\bf A}(p)$ such that $\langle V,W\rangle =C$ with $-1 < C < 1$ (one could take $C\geq 0$) and such that there exists sequences $V_i\in {\bf A}(p)$ and $W_i\in {\bf A}(p)$ such that $V_i$ converges to $V$ and $W_i$ converges to $W$.
We extend $V_i$ and $W_i$ to vector fields by letting $V_i(q)=\nabla b_{1V_i}(q)$ and $W_i(q)=\nabla b_{1W_i}(q)$.  By the uniform bound on Lipshitz constants on Busemann functions, i.e. Lemma \ref{busemannsmooth}, there is an $\epsilon>0$ (depending only on $C$ and the Lipshitz constant of the busemann functions but not on $i$) such that for all $q\in B_p(\epsilon)$ (the $\epsilon$ ball about $p$) and all sufficiently large $i$ we have
$$-\frac{1+C}{2} < \langle V_i(q),W_i(q)\rangle < \frac{1+C}{2}.$$
This holds since for large $i$ we have $\langle V_i(p),W_i(p)\rangle $ is approximately $C$ and then, with respect to a parallel frame along geodesics of length $\epsilon$, the change of $V_i$ and $W_i$ is uniformly bounded by the Lipshitz constant and $\epsilon$.  We can further take $\epsilon$ less than the distance from $p$ to $\partial \widetilde{M_1}$.

Now consider the Busemann function $b_{1V_i}$ on the $0$ level set $H_{W_i}(0)$ of $b_{1W_i}$ (i.e. the level set through $p$).  By the inner product condition above we can find unit speed differentiable curves $\tau_1$ and $\tau_2$ in $H_{W_i}(0)$ starting at $p$ of length $\epsilon$ (and hence in $H_{W_i}(0)\cap B_p(\epsilon)$) such that
$$\langle \tau_1'(s),V_i(\tau_1(s))\rangle > \bar C=\sqrt{1-\Big (\frac{1+C}{2}\Big )^2}\ $$
and
$$\langle \tau_2'(s),V_i(\tau_2(s))\rangle < -\bar C=-\sqrt{1-\Big (\frac{1+C}{2}\Big )^2}\ .$$
Thus for every sufficiently large $i$ there are points $z^1_i,z^2_i \in H_{W_i}(0)\cap B_p(\epsilon)$ such that $b_{1V_i}(z^1_i) > \bar C \epsilon$ and $b_{1V_i}(z^2_i) < -\bar C \epsilon$.  Thus by lemma \ref{levelsets} for every sufficiently large $i$ there are points $y^1_i$ and $y^2_i$ on the boundary of $H_{W_i}(0)\cap \widetilde{M_1}$ with $b_{1V_i}(y^1_i) > \bar C \epsilon$ and $b_{1V_i}(y^2_i) <-\bar C \epsilon$.

Since $V_i$ and $W_i$ converge to trapped geodesics, Lemma \ref{vertical} says that as $i\to \infty$ the lines $L_{V_i}$ and $L_{W_i}$ converge to vertical.  But this means that the the Busemann functions $b_{1V_i}$ and $b_{1W_i}$ (which are height functions outside $\widetilde M_1$) approximate the vertical height function.  In particular, for $i$ large enough the values of $b_{1V_i}$ on the boundary of $H_{W_i}(0)\cap \widetilde{M_1}$ vary by less than $\bar C\epsilon$.  This contradicts the simultaneous existence of $y^1_i$ and $y^2_i$ for large $i$.

\qed

The first consequence of this proposition is that the trapped geodesics are also minimizing (as limits of minimizing geodesics) and hence $g$ has no conjugate points.  Start with a large enough flat $n+1$ torus $T^n\times S^1$ so that $D^n$ sits isometrically in $T^n$.  Now if we replace $D^n\times S^1$ with $(M_1,\partial M_1,g)$ we get an $n+1$ torus since $M_1$ is diffeomorphic to $M$ (as was pointed out after the proof of Lemma \ref{geosminimize}).  Further it has no conjugate points.  Then the theorem of Burago-Ivanov \cite{Bu-Iv94} proving the E. Hopf conjecture says that the metric is flat.  This gives a proof of Theorem \ref{3+dims}.

However the above proof does not generalize very far.  In the next section we give an alternative proof that does generalize.

\section{generalizations}
\label{gens}

In the previous section we considered only flat metrics so as to make the the proof more transparent.  However the arguments extend almost without change to give

\begin{proposition}
\label{mainhighergen}
 Let $D^n$ be a ball in a complete simply connected manifold $N^n$ with nonpositive curvature, and $(M_1,\partial M_1,g_1)$ a Riemannian manifold with boundary that has the same scattering data as $(D\times S^1,\partial D\times S^1,g)$ where $g$ is the product metric.  Then through every point of $M_1$ there is exactly one trapped geodesic.
\end{proposition}

The main change that affects the proof is that geodesics in $\widetilde{M_1^{ext}}$  are not lines outside $\widetilde M_1$ but geodesics in $N$.  Oriented geodesic rays in $N$ are thought of as parallel if they have the same limit point at infinity (which means that they stay a bounded distance from each other).  The notion of ``positive'' geodesics also makes sense.  This allows us to relate ${\bf A}(p)$ to ${\bf A}(q)$.  Most of the arguments go through exactly as before.  In particular Lemmas \ref{fundgroup}, \ref{lensdata}, and \ref{geosminimize} go through as is.  Lemma \ref{vertical} also goes through as before where convergence of directions needs to be interpreted as the endpoints at infinity converging (where the topology of infinity is the cone topology - hence homeomorphic to the standard $n-1$ sphere).   The argument in Lemma \ref{levelsets} needs to be viewed carefully since there may be many oriented geodesics corresponding to $-\nabla b_{1W}$ as the point on $H_W(s)$ varies.  However, none of these will be positive so again the proof goes through.  The only part of Proposition \ref{mainhigher} that needs to be noted is that in nonpositive curvature the level sets of Busemann functions vary continuously with the initial vector.

\qed

\begin{remark}
The above arguments likely can be modified to cover manifolds without conjugate points.   One first needs to deal with the fact that balls may not be convex.  This would seem to give us problems with the differentiability of the metric on $\widetilde{M_1^{ext}}$ at the boundary of $\widetilde M_1$.  However (except possibly in the case where the boundary of $D$ contains a region that is totally geodesic) since there are no conjugate points Theorem 1 of \cite{St-Uh09} will still tell us that the metric will be smooth.  In fact very little of the argument really needs the boundary to be convex (or the metric to be smooth for that mater) but extending the arguments would look somewhat different from the above.  Also one has to worry about relating ${\bf A}(p)$ and ${\bf A}(q)$.  This can be done by fixing a base point $x_0\in \widetilde{M_1^{ext}}-\widetilde M_1$ and looking only at Busemann functions defined by vectors in $U_{x_0}$.  Also Busemann functions are not as well behaved.  In any event, the arguments would look very different and we wont pursue them here.
\end{remark}

\bigskip

We now want to generalize Theorem \ref{3+dims} to the nonflat case.  The first point to note is that for $g$ and $g_1$ as in Proposition \ref{mainhighergen} we have
$$Vol(g_1)=Vol(g).$$

Let ${\bf T_1}\subset UM_1$ be the set of unit vectors tangent to trapped geodesic rays.  Similarly define ${\bf T}\subset UM$ (here $M=D\times S^1$).
Using the fact that the metrics are lens equivalent we consider the standard measure preserving map $F:UM_1-{\bf T}\to UM-{\bf T}$ which assigns to each vector $V\in UM_1-{\bf T}$ the unique vector $W\in UM_1-{\bf T}$ such that $V=\gamma'_1(t)$ and $W=\gamma'(t)$ where $\gamma$ and $\gamma_1$ are geodesics with corresponding initial conditions on the boundary.  I.e. $\gamma(0)\in \partial M$ and $\gamma_1(0) \in \partial M_1$ are corresponding points while $\gamma'(0)$ and $\gamma'_1(0)$ are corresponding inwardly pointing unit vectors.  That $F$ (which conjugates the geodesic flows) is measure preserving is a standard fact which follows for example from Santal\'o's formula (see for example \cite{Cr91}).  The fact that ${\bf T}$ and ${\bf T_1}$ have measure 0 tells us that the unit tangent bundles (and hence the manifolds) have the same volume.

The generalization of Theorem \ref{3+dims} is

\begin{theorem}
\label{generalizations}

Let $D^n$ be a ball in $N^{n-1}\times \R$ where $N$ is a complete simply connected manifold with nonpositive curvature.  Then $D\times S^1$ is scattering rigid.

\end{theorem}

{\bf Proof:}
One proves this via Proposition 2.2 of \cite{Cr-Kl98}.  That result compares two metrics $g$ and $g_1$ on manifolds without conjugate points with an additional condition on $g$ that it contain a real factor (which is satisfied by our assumption on $D$). Under a weak lens equivalency condition (which is satisfied since the metrics are lens equivalent) one concludes that $Vol(g_1)\geq Vol(g)$ with equality holding if and only if $g_1$ is isometric to $g$.  This (along with the above fact that $Vol(g_1)=Vol(g)$) proves the Theorem.

\qed

Note that we do not claim the the result for $D^n$ a ball in a complete simply connected manifold with nonpositive curvature.  $D^n$ is a ball in $N^{n-1}\times \R$.  This means that in the universal cover of $D$ (and hence of $M=D\times S^1$) there is a bounded $\R$ factor (as opposed to the unbounded $\R$ factor coming from the $S^1$ factor in $M$).   This is because Proposition 2.2 of \cite{Cr-Kl98} requires a set whose $\R$ component in the universal cover is bounded.

%\pagebreak

\end{document}